\documentclass[a4paper,12pt]{amsart}
\usepackage[utf8]{inputenc}
\usepackage{amsmath,amssymb,amsfonts,amsthm,graphicx}

\newtheorem{theorem}{Theorem}[section]

\newtheorem{corollary}[theorem]{Corollary}
\numberwithin{equation}{section}

\theoremstyle{remark}

\def\ppt{\frac{\partial}{\partial t}}
\def\pps{\frac{\partial}{\partial s}}

\title{Harnack inequalities for the curve shortening flow}

\author{Mihai B\u{a}ile\c{s}teanu}
\thanks{Department of Mathematical Sciences, Central Connecticut State University, 120 Marcus White Hall, New Britain, CT 06050, USA \texttt{mihaib@ccsu.edu}}

\begin{document}

\begin{abstract}
We offer an algorithmic approach for determining Harnack quantities for the curve shortening flow and we show how, following this procedure, one can obtain Hamilton's Harnack inequality for this flow $\kappa_t+\frac{1}{2t}\kappa\geq\frac{\kappa_s^2}{\kappa}$, where $\kappa$ is the curvature of the curve being deformed by the flow.   
\end{abstract}

\maketitle

\section{Introduction}

Given a family of planar curves $\gamma(.,t)$ and a choice of continuous unit normal vectors $\mathbf{N}(.,t)$, the evolution equation \[\begin{cases}\frac{\partial \gamma}{\partial t}=\kappa \mathbf{N}\\
\gamma(.,0)=\gamma_0(.)
\end{cases}\]
defines the curve shortening flow ($\kappa$ denotes the curvature). In other words, this describes the motion of a curve when each point on the curve moves in the direction of the normal vector, with speed equal to the curvature at that point. The equation has a unique smooth solution (\cite{gage-hamilton86}) on some time interval $[0,T)$.  

Several important techniques which are typical for the field of geometric flows appear in the study of the curve shortening flow in a less technical and more intuitive way. This is because, for example, this flow is the 1-dimensional version of the mean curvature flow. 

Historically, the curve shortening flow was first proposed in 1956 by Mullins to model the motion of idealized grain boundaries (\cite{mullins56}). The study of this flow developed during the 1980's, through the work of M. Gage (see \cite{gage83,gage84,gage-epstein87,gage90,gage90-1,gage93,gage94}) and R. Hamilton (see \cite{gage-hamilton86,hamilton95}) on convex plane curves and M. Grayson (see \cite{grayson87,grayson89} on embedded plane curves, at the time when geometers started using geometric flows to study topological problems.   

In particular, Gage and Hamilton proved that convex curves shrink to a point in finite time and they become asymptotically round while approaching that point. Later, Grayson showed that, in fact, any closed embedded curve in the plane becomes convex eventually, and thus converges to a point. The solution is smooth, without singularities, in all these results. Later the proof has been simplified, using isoperimetric estimates by Huisken \cite{huisken98} and Hamilton \cite{hamilton92}. Moreover, Hamilton's Harnack estimate for the mean curvature flow in \cite{hamilton95} can be adapted to the curve shortening flow to prove Grayson's result. Recently, a completely new proof of Grayson's theorem, based also on isoperimetric estimates, was 
presented by Andrews and Bryan \cite{andrews-bryan11}. 

In this paper, we are giving a new proof of Hamilton's Harnack inequality involving the curvature $\kappa$. In \cite{hamilton95} Hamilton proved a Harnack estimate for the mean curvature flow, which when applied to the curve shortening flow (note that this is just the mean curvature flow in dimension one) becomes the following theorem: 

\begin{theorem}\label{thm-Harnack}
 For any closed convex curve in the plane changing under the curve shortening flow, one has that 
 \[\kappa_t+\frac{1}{2t}\kappa\geq\frac{\kappa_s^2}{\kappa}.\]
\end{theorem}

Here $s$ is the arc-length parameter of the curve and $\kappa$ is the curvature. 

The proof of this theorem is based on the maximum principle. One denotes the Harnack quantity $h(s,t)=\kappa_t+\frac{1}{2t}\kappa-\frac{\kappa_s^2}{\kappa}$ and after calculating the time evolution of $h$ (which is also a heat-type equation), one uses the maximum principle to show that $h\geq 0$ at all times on a closed curve. The proof does not show, however, how to actually come up with the Harnack expression. 

In the field of geometric flows, finding the precise Harnack quantity is usually the key in proving the Harnack estimates. Often times there isn't any algorithm which leads to the desired expression, it's a trial and error process. For example, Hamilton's way of obtaining Harnack guantities is to analyze solitons (self-similar solutions) of the flow and find expression involving the solitons which are zero.  

In this paper, we offer an algorithmic approach for determining Harnack quantities and we also show how, for a certain choice of constants, one can obtain precisely Hamilton's expression. The algorithm both determines the Harnack quantity and proves that it's non-negtive, using the maximum principle. The method is quite technical, but it has the advantage of being constructive.    

The method can be summarized as follows: let $f>0$ denote the function for which one needs to prove a Harnack inequality. Usually $f$ satisfies a heat-type equation. Denote with $u=\log f$, then $u_t=\triangle u+|\nabla u|^2$+other terms. Inspired by the expression of $u_t$, we define the Harnack quantity $H=a\triangle u+b |\nabla u|^2+\text{possibly other terms}+\varphi(x,t)$. $\varphi$ is a function that needs to be determined and that has to go to infinity as $t\to 0$. Next we assume that there is a time when $H<0$. Since $\varphi$ is very large at time $0$, it means that  $\lim\limits_{t\to 0^+} H=\infty$, i.e. $H$ is for sure positive at time $0$. There must be thus a time when $H=0$. At that time we can apply the maximum principle. In that sense, we determine $H_t-\triangle H -2\nabla u\nabla H$ (which has to be non-negative) and we impose conditions on $a, b$ and restrictions for $\varphi$ (from which one can build $\varphi$) such that the expression turns out positive, contradicting the initial assumption. Hence, in the same time, we both determine $H$ (by determining $a$, $b$ and building $\varphi$) and show that it's non-negative. 

This method proved to be successful in the study of Ricci flow and other heat-type equations (see, for example, the work of Cao \cite{cao-hamilton09,cao08,cao-ljun-liu13, cao-cerenzia-kazaras14}).

Our main result can be stated as follows:

\begin{theorem}\label{thm-princ}
Let $\gamma(.,t)$ be a family of closed convex curves being deformed under the curve shortening flow. Let $\kappa(s,t)$ denote the curvature of the curve at point $s$ and time $t$ and let $u(s,t)=\log\kappa(s,t)$. Then the expression
\[h_{\epsilon}(s,t):=u_{ss}+e^{2u}+\left(\frac{1}{2}+\epsilon\right)\frac{1}{2}\]
satisfies $h_{\epsilon}(s,t)\geq 0$ for any time $t>0$ and any $\epsilon>0$. 
\end{theorem}

Note that the $h_{\epsilon}\geq 0$ is actually equivalent to $\kappa_t+\left(\frac{1}{2}+\epsilon\right)\frac{\kappa}{t}\geq\frac{\kappa^2_s}{\kappa}$, so by taking $\epsilon\to 0$, we obtain, as as a corollary, Hamilton's result. 

The paper is structured as follows: in section 2 we present the setting of the problem and various results related to the curve shortening flow that will be used in our proof. Section 3 shows the above algorithm applied to this setting and as a result, the expression of the Harnack quantity is determined. For completeness of the exposition, we prove the non-negativity of the expression in section 4 and in section 5 we show how it relates to Hamilton's inequality.

\textbf{Acknowledgements} The author would like to thank Prof. Xiaodong Cao for suggesting to apply this method to the curve shortening flow and to Melanie Stam for useful discussions.   

\section{Background and setting}

We consider a family of parametrized immersed curves $\gamma(x,t)$, with $t>0$. This means that, for each fixed $t$, $\gamma: I\to\mathbb{R}^2$ is a smooth map from an interval $I\subseteq\mathbb{R}$ satisfying the condition: $|\gamma'(x)| \neq 0$ for every $x\in I$. Moreover the curves are invariant under both reparametrization of the domain $\mathbb{R}$ and rigid motions of the range $\mathbb{R}^2$.

The curves are considered to be embedded (i.e. they are homeomorphisms to their images). Intuitively, it's a curve that doesn't intersect itself. Another condition we impose on the curves is that they are closed (the maps are periodic with some period).

We parametrize the curves by arc length (the length of the curve becomes 1 under this parametrization) and we call the arc length parameter $s$. The unit tangent vector $\mathbf{T}$ is defined as 
\[\mathbf{T}=\frac{d\gamma}{ds}=\frac{\gamma'}{|\gamma'|}\]
where $|.|$ is the length of the vector. 

The vector function $\mathbf{T}$ has length one, hence its derivative will be perpendicular to it, and thus is a multiple of the unit normal vector $\mathbf{N}$. This gives the curvature and the Frenet-Serret formulas for planar curves:
\begin{align*}
\frac{d\mathbf{T}}{ds} & =\kappa \mathbf{N}\\
\frac{d\mathbf{N}}{ds} & =-\kappa \mathbf{T}
\end{align*}
where $\kappa$ is the curvature at that particular point. Note that for a convex curve $\kappa\geq 0$ at every point. 

If one wants to deform curves in the plane, a natural choice would be the curve shortening flow, which is a heat type equation defined as:
\begin{align}\label{CSF-eqn}
\frac{\partial \gamma(s,t)}{\partial t}=\kappa(s,t) \mathbf{N}(s,t).
\end{align}
with an intial condition $\gamma(s,0):=\gamma_0(s)$. To simplify the notation, we will drop the arguments in the future. 

It is not obvious that this is a heat type equation, but this fact becomes clear when noting that (\ref{CSF-eqn}) is equivalent to 
\[\frac{\partial \gamma}{\partial t}=\frac{\partial^2\gamma}{\partial s^2}.\]
But since the arc length parameter is defined using $\gamma$, the equation is not, in fact, linear ($s$ depends on $t$ also). 

Under this flow, the time-derivative does not commute with the derivative in $s$, in fact (cf. \cite{gage-hamilton86}, Lemma 3.1.3):
\[\frac{\partial}{\partial t}\frac{\partial}{\partial s}=\frac{\partial}{\partial s}\frac{\partial}{\partial t}+\kappa^2\frac{\partial}{\partial s}\]

As a result of this, the curvature evolves under the curve shortening flow according to the equation (\cite{gage-hamilton86}, Lemma 3.1.6)
\begin{align}\label{curv-eqn}
\frac{\partial \kappa}{\partial t}=\frac{\partial^2\kappa}{\partial s^2}+\kappa^3.  
\end{align}

The significant consequence of this last equation is that if the curve is convex at the beginning of the flow ($\kappa>0$ at time $t=0$), then it will stay convex as long as the flow exists.  

From now on, we will denote $\frac{\partial \kappa}{\partial t}=\kappa_t$, $\frac{\partial \kappa}{\partial s}=\kappa_s$, $\frac{\partial}{\partial s}=\partial_s$ and so on. 

We also assume that the initial curve is convex.

\section{Finding the Harnack quantity}

We start with the equation \ref{curv-eqn} satisfied by the curvature $\kappa$:
\begin{align*}
\kappa_t=\kappa_{ss}+\kappa^3 
\end{align*}
We assume that $\kappa>0$ at time 0, which implies that the curve will remain convex at all times. This can be shown by applying the maximum principle to equation (\ref{curv-eqn}). Hence $\kappa> 0$ at all times. 

This fact allows us to define $u(s,t)=\log \kappa(s,t)$. Notice that $u$ satisfies the following equation: \[u_t=u_{ss}+u_s^2+e^{2u}\]

Next we consider the following Harnack quantity:
\[h(s,t)=a u_{ss}+b u_s^2+ce^{2u}+\varphi(s,t)\]
where $a,b,c\in\mathbb{R}$. These numbers will be determined later, together with the function $\varphi$. 

We will compute the time evolution of $h$ and then we will use the maximum principle to establish the positivity of $h$ given a particular choice of the function $\varphi(s,t)$. The function $\varphi(s,t)$ has to be very large at $t=0$, in order to dominate the other terms and ensure that at time close to $0$ the Harnack quantity is positive. 

To find $h_t$ we begin by computing $(u_{ss})_t$ and $(u_s^2)_t$. Recall that the vector fields $\ppt$ and $\pps$ don't commute, in fact $[\partial_t, \partial_s]=k^2 \partial_s$.
\[(u_{ss})_t=\partial_t \partial_s \partial_s u= \partial_s \partial_t \partial_s u +k^2 u_{ss}=\partial_s\partial_s\partial_t u +\partial_s(k^2 u_s )+k^2 u_{ss}= (u_t)_{ss}+2k k_s u_s +2 k^2 u_{ss}\]

Further, one can replace $u_t$ with $u_{ss}+u_s^2+e^{2u}$ and obtain:
\begin{align*} 
(u_t)_{ss} & =\frac{\partial^2}{\partial s^2}(u_{ss}+u_s^2+e^{2u})\\
           & = u_{ssss}+(u_s^2)_{ss}+2e^{2u}u_{ss}+4e^{2u}u_{s}^2 
\end{align*}
Therefore: 
\begin{align}\label{u-sss-t}
 (u_{ss})_t =u_{ssss}+ (u_s^2)_{ss} +4e^{2u}u_{ss}+6e^{2u}u_s^2 
\end{align}

Next, we calculate the evolution of $u_s^2$:
\[(u_s^2)_t=2u_s u_{st}=2u_s(u_t)_s+2k^2u_s^2=2u_s(u_t)_s+2e^{2u}u_s^2\]
Using again the relationship $u_t=u_{ss}+u_s^2+e^{2u}$ leads to:
\begin{align*}
2 u_s \cdot (u_t)_s & = 2u_s\cdot\partial_s(u_{ss}+u_s^2+e^{2u})=2u_su_{sss}+4u_s^2u_{ss}+4e^{2u}u_s^2 \\
     & = (u_s^2)_{ss} -2u_{ss}^2 +4u_s^2u_{ss}+4e^{2u}u_s^2
\end{align*}
It follows that:
\begin{align}\label{u-s-2-t}
 (u_s^2)_t & =(u_{s}^2)_{ss} -2u_{ss}^2 +2(u_s^2)_s\cdot u_s +6e^{2u}u_s^2 \\ \nonumber   &=(u_{s}^2)_{ss} -2u_{ss}^2 +4 u_s^2u_{ss} +6e^{2u}u_s^2 
\end{align}

We are now ready to determine the evolution of the Harnack quantity. Using the expressions in (\ref{u-sss-t}) and (\ref{u-s-2-t}) leads to:

\begin{align*} 
h_t & =a(u_{ss})_t+b(u_s^2)_t+2ce^{2u}u_t+\varphi_t(s,t) \\
    & =a(u_{ssss}+ (u_s^2)_{ss} +4e^{2u}u_{ss}+6e^{2u}u_s^2)+b((u_{s}^2)_{ss} -2u_{ss}^2 +4 u_s^2u_{ss} +6e^{2u}u_s^2)\\
    &+2ce^{2u}(u_{ss}+u_s^2+e^{2u})+\varphi_t(s,t)\\
    & = \left[a u_{ssss}+b (u_s^2)_{ss}+4ce^{2u}u_s^2+2ce^{2u}u_{ss}+\varphi_{ss}\right] \\
    &+2(a u_{sss}+2b u_s u_{ss}+2ce^{2u}u_s+\varphi_s)u_s\\
    &-6ce^{2u}u_s^2-\varphi_{ss}+2a u_{ss}^2-2b u_{ss}^2+6b e^{2u}u_s^2-2\varphi_su_s+4a e^{2u}u_{ss}+6a e^{2u}u_s^2+2ce^{4u}+\varphi_t\\
    &=h_{ss}+2h_su_s+4e^{2u}h+2(a-b)u_{ss}^2+(6a+2b-6c)e^{2u}u_s^2-2ce^{4u}-4\varphi e^{2u}-\varphi_{ss}+\varphi_t-2\varphi_su_s    
\end{align*}

Therefore the Harnack quantity $h$ evolves according to this expression:
\begin{align}\label{h-t}
h_t & =h_{ss}+2h_su_s+4e^{2u}h+2(a-b)u_{ss}^2\\
\nonumber & +(6a+2b-6c)e^{2u}u_s^2-2ce^{4u}-4\varphi e^{2u}-\varphi_{ss}+\varphi_t-2\varphi_su_s
\end{align}

Assume there is a time when $h<0$. Since $\varphi(s,t)$ is assumed to go to infinity at time $0$, $h$ has to be positive close to the start of the flow. But since the solution is smooth (it's a heat-type equation), $h$ is also smooth, so there has to be a  first time $t_0$ when $h(s_0,t_0)=0$. It follows from the maximum principle that $h_t\leq 0$, $h_s=0$ and $h_{ss}\geq 0$. Moreover, assuming that $a\neq 0$, at this point
\[u_{ss}=-\frac{b u_s^2+ce^{2u}+\varphi}{a}.\]

As a consequence, at $t_0$ the expression (\ref{h-t}) becomes:
\begin{align}\label{ineg-init}
0\geq & \frac{2(a-b)}{a^2}[b u_s^2+ce^{2u}+\varphi]^2 \\
      & + (6a+2b-6c)e^{2u}u_s^2-2ce^{4u}-4\varphi e^{2u}-\varphi_{ss}+\varphi_t-2\varphi_su_s \nonumber
\end{align}

To simplify the calculation, denote $e^{2u}=X>0$ and $u_s^2=Y\geq 0$, so the above can be rewritten as:
\begin{align*}
0\geq & \frac{2(a-b)}{a^2}[b Y+cX+\varphi]^2 \\
      & + (6a+2b-6c)XY-2cX^2-4\varphi X-\varphi_{ss}+\varphi_t-2\varphi_su_s\\
      &= \frac{2(a-b)}{a^2}\left[b^2 Y^2+\left(c^2-\frac{ca^2}{a-b}\right)X^2\right]+\left[6a+2b-6c+\frac{4(a-b)b c}{a^2}\right]XY\\
      &+\left(\frac{4c(a-b)}{a^2}-4\right)\varphi X+\frac{4b(a-b)}{a^2}\varphi Y-\varphi_{ss}+\varphi_t-2\varphi_su_s+\frac{2(a-b)}{a^2}\varphi^2
\end{align*}

We now proceed to analyze each term and impose conditions on $a,b,c$ and $\varphi$ in order to make the expression on the right positive (the goal is, in fact, to produce a contradiction). We will assume each term is non-negative, and at the end, after we have analyzed all the conditions, we will set one value to be positive, thus making the whole expression positive. 

\bigskip
1) For the first term, let us start by assuming $a>b\geq0$. This will assure that $2(a-b)$ is positive. Next, we need \[c^2-\frac{ca^2}{a-b}\geq 0\]
This is equivalent to $c\left(c-\frac{a^2}{a-b}\right)\geq 0$, so either $c\leq 0$ or $c\geq\frac{a^2}{a-b}>0$.

\bigskip
2) Since both $X,Y$ are non-negative, the non-negativity of the second term is satisfied provided that \[6a+2b-6c+\frac{4(a-b)b c}{a^2} \geq 0\] which is equivalent to 
\[c\leq \frac{(3a+b)a^2}{3a^2-2ab+2b^2}\]
Notice that both the numerator and the denominator are positive ($3a^2-2ab+2b^2=2a^2+b^2+(a-b)^2$), therefore if $c\leq0$ we get this condition for free.

\bigskip
3) For the third term, as $X$ is positive and $\varphi$ will be chosen to be positive, we need $\frac{4c(a-b)}{a^2}-4\geq 0$. This gives $c\geq \frac{a^2}{a-b}$, which is consistent with $1)$. Therefore from the analysis of the first three terms we obtain
\begin{align}\label{c-exp}
 \frac{a^2}{a-b}\leq c\leq \frac{(3a+b)a^2}{3a^2-2ab+2b^2}
\end{align}

\bigskip
4) We group the last four terms into one expression 
\begin{align}\label{phi-init}
\frac{4b(a-b)}{a^2}\varphi Y-\varphi_{ss}+\varphi_t-2\varphi_su_s+\frac{2(a-b)}{a^2}\varphi^2 
\end{align}

We will find a positive function $\varphi(s,t)$ of two variables, which would turn the above expression positive. 

Proceeding as in, for example, \cite{cao-cerenzia-kazaras14}, by Cauchy-Schwartz one obtains that:
\[\frac{4b(a-b)}{a^2}\varphi Y-2\varphi_su_s=\frac{4b(a-b)}{a^2}\varphi u_s^2-2\varphi_su_s\geq -\frac{a^2\varphi_s^2}{4b(a-b)\varphi}\]
effectively getting rid of the term $-2\varphi_su_s$, which is harder to deal with on its own. 

Therefore, the expression (\ref{phi-init}) is greater than:
\begin{align}\label{phi-exp}
\varphi_t-\varphi_{ss}-\frac{a^2\varphi_s^2}{4b(a-b)\varphi}+\frac{2(a-b)}{a^2}\varphi^2 
\end{align}

There is one more condition that needs to be imposed on $\varphi$: as $(t,s)\to (0,0)$, $\varphi\to\infty$. This assures that the Harnack quantity $h(s,t)$ is very large at time $0$ and so it is positive at the beginning of the flow. 

A reasonable ansatz is $\varphi(s,t)=\frac{\alpha}{t}+\frac{\beta}{s^2}$ with $\alpha\geq0$, $\beta\geq0$ (to be determined later). This means that:
\begin{align*} 
\varphi &= \frac{\alpha}{t}+\frac{\beta}{s^2} & \varphi_t &= -\frac{\alpha}{t^2}\\
\varphi_s &= -\frac{2\beta}{s^3} & \varphi_{ss} &= \frac{6\beta}{s^4}
\end{align*}
 
Notice that $\varphi\geq \frac{\beta}{s^2}$, hence 
\[\frac{\varphi_s^2}{\varphi}=\frac{4\beta^2}{\varphi s^6}\leq \frac{4\beta}{s^4}.\]

Denoting $A=\frac{2(a-b)}{a^2}$ and $B=\frac{a^2}{4b(a-b)}$, (\ref{phi-exp}) becomes:
\begin{align}\label{varphi-pos}
\varphi_t-\varphi_{ss}-B \frac{\varphi_s^2}{\varphi}+A\varphi^2 
\end{align}

We plug the expressions of $\varphi_s$, $\varphi_{ss}$ and $\varphi_t$ into (\ref{varphi-pos}): 
\begin{align}\label{alf-bet}
 \varphi_t-\varphi_{ss}-B \frac{\varphi_s^2}{\varphi}+A\varphi^2 &= -\frac{\alpha}{t^2} - \frac{6\beta}{s^4}-B \frac{\varphi_s^2}{\varphi}+ A\left( \frac{\alpha}{t}+\frac{\beta}{s^2}\right)^2\nonumber\\
                  & \geq  -\frac{\alpha}{t^2} - \frac{6\beta}{s^4}- \frac{4\beta B}{s^4}+ A\left( \frac{\alpha}{t}+\frac{\beta}{s^2}\right)^2\nonumber\\
                  & =\frac{A\alpha^2-\alpha}{t^2}+\frac{2A\alpha \beta}{ts^2}+\frac{A\beta^2-6\beta-4\beta B}{s^4}
\end{align}

By choosing $\alpha$ and $\beta$ large enough (in particular $\alpha\geq\frac{1}{A}$ and $\beta\geq\frac{6+4B}{A}$ and at least one of them being strict), all the terms in the above expression are positive, which would give that 
\[0< \varphi_t-\varphi_{ss}-\frac{a^2\varphi_s^2}{4b(a-b)\varphi}+\frac{2(a-b)}{a^2}\varphi^2\]

contradicting the expression in (\ref{ineg-init}). This means that there isn't any point where $h(s,t)<0$ and since $\lim\limits_{t\to 0}\varphi(s,t)=\infty$, it means that $h(s,t)$ always stays non-negative. 

Let's summarize all the conditions on $a,b,c$ and $\varphi$:
\begin{itemize}
 \item $a,b,c\geq 0$
 \item $a>b\geq 0$
 \item $\frac{a^2}{a-b}\leq c\leq \frac{(3a+b)a^2}{3a^2-2ab+2b^2}$
 \item $\varphi(s,t)=\frac{\alpha}{t}+\frac{\beta}{s^2}$
 \item $\alpha\geq\frac{a^2}{2(a-b)}$
 \item $\beta\geq \frac{a^2(a^2+6b(a-b))}{2b(a-b)^2}$
\end{itemize}

Focusing on the restrictions on $c$, one will notice that \[\frac{a^2}{a-b}\leq \frac{(3a+b)a^2}{3a^2-2ab+2b^2}\] which is equivalent to
\[3a^4-2a^3b+2a^2b^2\leq 3a^4+a^3b-3a^3b-a^2b^2\] or
\[3a^2b^2\leq 0\]
The only possibility for this to be true is if $b=0$ (recall that $a>b\geq0$), which forces $c=a$ and $\alpha\geq \frac{a}{2}$. For $\beta$ one has to choose $\beta=0$, otherwise it would be undefined (notice that this actually makes sense, since the expression in (\ref{alf-bet}) still remains positive provided that $\alpha$ is slightly larger than $a/2$, i.e. $\alpha=a/2+\epsilon$).

Putting all these together allows us to conclude that a reasonable choice for the Harnack quantity following this procedure is:
\[h(s,t)=a u_{ss}+a e^{2u}+\left(\frac{a}{2}+\epsilon\right)\frac{1}{t}\]
for some positive $a$ and $\epsilon$. We can rescale the quantity and pick $a=1$, therefore the Harnack quantity is \[h(s,t)=u_{ss}+e^{2u}+\left(\frac{1}{2}+\epsilon\right)\frac{1}{t}\] for an arbitrarily small $\epsilon>0$.

\textbf{Discussion about the choice of $\varphi(s,t)$} Notice that the choice of $\varphi(s,t)$ is not unique. We made the ansatz $\varphi(s,t)=\frac{\alpha}{t}+\frac{\beta}{s^2}$ because it's the simplest function which goes to $+\infty$ as $(s,t)$ approaches $(0,0)$. If the reader wonders why $s^2$ was chosen, instead of $s$, the answer is that this came from the fact that we needed to estimate $\frac{\varphi_s^2}{\varphi}$ and relate it to $\varphi_{ss}$. It is by no means the only solution. In fact, it's worth investigating what kind of inequalities one gets by changing the expression of $\varphi$. 

\section{The proof}

With the choice of $h(s,t)$ from the previous section, we will prove theorem \ref{thm-princ}. This is done for completeness of exposition and as a check-up for the expression found following the procedure. 

Let $h(s,t):=u_{ss}+e^{2u}+\left(\frac{1}{2}+\epsilon\right)\frac{1}{t}$ be the Harnack quantity. Assume there is a point $(s_1,t_1)$ such that $h(s_1,t_1)<0$. Since $\lim\limits_{t\to 0}h(s,t)=\infty$ and $h$ is smooth, then there has to be a first point $(s_0,t_0)$ where $h$ is actually 0. 

Computing the time evolution of the Harnack quantity at $t_0$, we obtain that:

\[h_t=h_{ss}+2h_su_s+4e^{2u}h+2u^2_{ss}-2e^{4u}-\frac{2}{t}e^{2u}-\frac{1}{2t^2}\]

We can now apply the maximum principle to the above expression: at the point $(s_0,t_0)$ the following holds: $h(s_0,t_0)=0$,  $h_t\leq 0$, $h_s=0$ and $h_{ss}\geq 0$. Therefore, at the point $(s_0,t_0)$:

\[0\geq 2u^2_{ss}-2e^{4u}-4\left(\frac{1}{2}+\epsilon\right)\frac{1}{t}e^{2u}-\left(\frac{1}{2}+\epsilon\right)\frac{1}{t^2}\]

Since at $(s_0,t_0)$, $h=0$, it means that $u_{ss}=-e^{2u}-\left(\frac{1}{2}+\epsilon\right)\frac{1}{t}$, therefore the above inequality becomes:

\begin{align*}
0 & \geq 2 \left[e^{2u}+\left(\frac{1}{2}+\epsilon\right)\frac{1}{t}\right]^2-2e^{4u}-4\left(\frac{1}{2}+\epsilon\right)\frac{1}{t}e^{2u}-\left(\frac{1}{2}+\epsilon\right)\frac{1}{t^2}\\
 & = \left[2\left(\frac{1}{2}+\epsilon\right)^2-\left(\frac{1}{2}+\epsilon\right)\right]\frac{1}{t^2}=(2\epsilon^2+\epsilon)\frac{1}{t^2}
\end{align*}
which is clearly a contradiction, since $\epsilon>0$ and $t>0$. 

Therefore there does not exist any time when $h$ is negative, hence $h(s,t)\geq0$ for all times $t>0$.

\section{Hamilton's Harnack inequality}

From theorem \ref{thm-princ}, we have that \[u_{ss}+e^{2u}+\left(\frac{1}{2}+\epsilon\right)\frac{1}{t}\geq 0\] where $u=\log\kappa$

By replacing $u_{ss}$ with $u_t-u_s^2-e^{2u}$, this is equivalent to:
\begin{align*}
u_t-u_s^2-e^{2u}+e^{2u}+\left(\frac{1}{2}+\epsilon\right)\frac{1}{t}\geq 0\\
\frac{\kappa_t}{k}-\frac{\kappa_s^2}{k^2}+\left(\frac{1}{2}+\epsilon\right)\frac{1}{t}\geq 0\geq 0\\
\kappa_t+\left(\frac{1}{2}+\epsilon\right)\frac{\kappa}{t}\geq \frac{\kappa_s^2}{\kappa}
\end{align*}

This is true for any $\epsilon>0$, so by taking the limit $\epsilon\to 0$, we obtain as a corollary, Hamilton's Harnack inequality:
\begin{corollary}
The curvature of any closed convex curve in the plane changing under the curve shortening flow satisfies 
 \[\kappa_t+\frac{1}{2t}\kappa\geq\frac{\kappa_s^2}{\kappa}.\]
\end{corollary}

\bibliographystyle{abbrv}
\bibliography{bio}

\end{document}